\def\wlog#1{}
 \makeatletter
 \renewcommand{\@latex@info}[1]{}
 \makeatother
 \renewcommand{\PackageInfo}[2]{}

 \makeatletter
 \renewcommand{\@font@info}[1]{}
 \makeatother

\documentclass[leqno, mathsec, numreferences, 12pt]{kluwer}

\usepackage{amscd}

\newcommand{\conservepaper}{
 \hoffset=-0.75in
 \setlength{\textwidth}{6.5in}
 \voffset=-0.5in
 \setlength{\textheight}{9.0in} %% only 9 inches for US paper
 }
 \conservepaper %% ???

\renewcommand{\copyrightline}{Submitted to \emph{Geometriae Dedicata}
 (May 2002).
 This version: 30 May 2002.} %% ???

\renewcommand{\thesubsection}{\thesection\Alph{subsection}}

\newcommand{\oldsection}{}
 \let\oldsection=\section
\renewcommand{\section}{\setcounter{equation}{0}\oldsection}
\renewcommand{\theequation}{\thesection.\arabic{equation}}

\newcommand{\operatorname}[1]{\mathop{\rm{#1}}}

\newcommand{\integer}{{\mathbb Z}}
\newcommand{\real}{{\mathbb R}}
\newcommand{\torus}{{\mathbb T}}
\newcommand{\rational}{{\mathbb Q}}

\newcommand{\Lie}[1]{\mathfrak{#1}}
\newcommand{\SL}{\operatorname{SL}}
\newcommand{\SO}{\operatorname{SO}}
\newcommand{\Ad}{\operatorname{Ad}}
\newcommand{\Hom}{\operatorname{Hom}}
\newcommand{\iso}{\cong}
\newcommand{\atimes}{\mathbin{\times_{\alpha}}}

\newcommand{\TPvert}{T(P)_{\normalfont{vert}}}

\renewcommand{\Lie}[1]{\mathfrak{\lowercase{#1}}}

\newcommand{\Rrank}{\mathop{\mbox{\normalfont$\real$-rank}}}

\newcommand{\eqref}[1]{{\upshape(\ref{#1})}}
\newcommand{\pref}[1]{{\upshape(\ref{#1})}}
\newcommand{\fullref}[2]{{\upshape\ref{#1}{\pref{#1-#2}}}}

\newtheorem{prop}[equation]{Proposition}
\newtheorem{thm}[equation]{Theorem}

\newtheorem{cor}[equation]{Corollary}

\newdisplay{defn}[equation]{Definition}
\newdisplay{notation}[equation]{Notation}
\newdisplay{rem}[equation]{Remark}
\newdisplay{eg}[equation]{Example}
\newdisplay{problem}[equation]{Problem}
\newdisplay{question}[equation]{Question}

\newproof{proof}{Proof}

 \newcounter{step}
\newenvironment{step}[1][\unskip]{\refstepcounter{step}
\em
 \medskip \noindent Step \thestep\
 #1.\ }{\unskip\upshape}

 \newcounter{case}

\newenvironment{thmref}{\thmrefer}{}
\newcommand{\thmrefer}[1]{\renewcommand\theequation
 {\protect\ref{#1}$'$}\addtocounter{equation}{-1}}

\newenvironment{fullthmref}{\fullthmrefer}{}
\newcommand{\fullthmrefer}[2]{\renewcommand\theequation
 {\protect\fullref{#1}{#2}$'$}\addtocounter{equation}{-1}}

\newcommand{\see}[1]{(see~\ref{#1})}
\newcommand{\cf}[1]{(cf.~\ref{#1})}

\begin{document}
 \begin{article}

\begin{opening}

\title{Ergodic actions of semisimple Lie groups 
 on compact principal bundles}

\author{Dave \surname{Witte}}

\institute{Department of Mathematics, Oklahoma State University,
Stillwater, OK 74078, USA
 \email{dwitte@math.okstate.edu,
  http://www.math.okstate.edu/$\sim$dwitte}
 }

\author{Robert J.~\surname{Zimmer}}
 \institute{Department of Mathematics, University of Chicago, Chicago, IL
60637, USA
 \email{r-zimmer@uchicago.edu} \\
Current address: Office of the
Provost, Brown University, Box 1862, Providence, RI 02912, USA
 \email{Robert\underline{ }Zimmer@brown.edu}}

\runningtitle{Ergodic actions on principal bundles}

\begin{abstract}
 Let $G = \SL(n,\real)$ (or, more generally, let $G$ be a connected,
noncompact, simple Lie group).  For any compact Lie group~$K$, it is easy
to find a compact manifold~$M$, such that there is a
volume-preserving, connection-preserving, ergodic action of~$G$ on some
smooth, principal $K$-bundle~$P$ over~$M$.  Can $M$ can be chosen
independent of~$K$?  We show that if $M = H/\Lambda$ is a homogeneous
space, and the action of~$G$ on~$M$ is by translations, then $P$~must
also be a homogeneous space $H'/\Lambda'$.  Consequently, there is a
strong restriction on the groups~$K$ that can arise over this
particular~$M$.
 \end{abstract}

\keywords{ergodic action, principal bundle, Borel cocycle, semisimple Lie
group, invariant connection}

\end{opening}

\section{Introduction}

The Margulis Superrigidity Theorem \cite[Thms.~7.6.6i, pp.~245--246, and
9.6.15i(a), p.~332]{MargulisBook} describes the finite-dimensional
representations of irreducible lattices in semisimple Lie groups of
higher real rank. In geometric terms, this is an explicit description of
a certain class of equivariant principal bundles. The following weak
version of this result records merely the conclusion that there are only
finitely many possible structure groups.

\begin{defn}
 A $G$-space~$M$ is \emph{irreducible} if every closed, noncompact,
normal subgroup of~$G$ acts ergodically on~$M$.
 \end{defn}

\begin{thm}[Margulis Superrigidity Theorem] \label{MargulisBundle}
 Given any connected, semisimple, linear Lie group~$G$, with 
 $\Rrank G \ge 2$,
 there is a corresponding finite set $\mathcal{H} = \{H_1,\ldots,H_n\}$
of connected, linear Lie groups, such that if
 \begin{itemize}
 \item $M = \Gamma \backslash G$ is an irreducible homogeneous $G$-space
of finite volume,
 \item $H$ is a connected, linear Lie group,
 \item $P$ is a principal $H$-bundle over~$M$,
 and
 \item the action of~$G$ on~$M$ lifts to an ergodic action of~$G$ on~$P$
by bundle automorphisms,
 \end{itemize}
 then $H$ is isomorphic to one of the groups in~$\mathcal{H}$.
 \end{thm}

Theorem~\ref{MargulisBundle} has been generalized to allow $G$-spaces
that are not homogeneous, under the additional assumption that $H$ is
semisimple, with no compact factors \cite[Thm.~5.2.5, p.~98]{ZimmerBook}.
(We note that $G$, being semisimple, has only finitely many normal
subgroups.)

\begin{thm}[Zimmer] \label{ZimmerBundle}
 If
 \begin{itemize}
 \item  $G$ is a connected, semisimple, linear Lie group, with 
 $\Rrank G \ge 2$,
 \item $M$ is an irreducible $G$-space with $G$-invariant finite volume,
 \item $H$ is a connected, semisimple, linear Lie group, with no compact
factors,
 \item $P$ is a principal $H$-bundle over~$M$,
 and
 \item the action of~$G$ on~$M$ lifts to an ergodic action of~$G$ on~$P$
by bundle automorphisms,
 \end{itemize}
 then $H$ is locally isomorphic to a normal subgroup of~$G$.
 \end{thm}

This latter result does not address the case where the semisimple Lie
group~$H$ is compact. It is perhaps surprising that the answer is
completely different in this situation. Indeed, there is no restriction
at all on the possible structure groups: for any compact Lie group~$K$,
it is easy to construct a $G$-equivariant principal $K$-bundle~$P$, on
which~$G$ acts irreducibly.

\begin{fullthmref}{unifK}{mfld}
 \begin{prop}
 Let $G$ be a noncompact, linear Lie group. For any compact Lie
group~$K$, there is a smooth, irreducible, volume-preserving action of~$G$
on some principal $K$-bundle~$P$ over some manifold~$M$.
 \end{prop}
 \end{fullthmref}

For our construction in \fullref{unifK}{mfld}, the manifold~$M$ depends
on the compact group~$K$. This suggests the following problem:

\begin{problem} \label{ProbTop}
 Let 
 \begin{itemize}
 \item $G$ be a connected, semisimple, linear Lie group, with no
compact factors, such that $\Rrank G \ge 2$, 
 and
 \item $M$ be a smooth manifold of finite volume, on which $G$~acts
irreducibly, by volume-preserving diffeomorphisms.
 \end{itemize}
 Find all the connected, compact Lie groups~$K$, for which there is a
principal $K$-bundle~$P$ over~$M$, such that the action of~$G$ on~$M$
lifts to an ergodic action of~$G$ on~$P$, by bundle automorphisms. In
particular, is there a choice of~$M$ for which every such group~$K$ is
possible?
 \end{problem}

Our construction in \fullref{unifK}{mfld} yields a principal bundle~$P$
with a $G$-invariant connection \cf{StandardEg}, so it is also natural to
consider the following geometric version of the problem:

\begin{problem}
 Consider Prob.~\ref{ProbTop}, with the additional requirement that
there is a $G$-invariant connection on the principal bundle~$P$.
 \end{problem}

We investigate this geometric question in the special case where $M$ is
one of the known actions of~$G$ that arise from an algebraic
construction. (It has been conjectured that every irreducible $G$-space
is isomorphic to one of these known actions, modulo a nowhere-dense,
$G$-invariant set. This may suggest that our results could be helpful in
understanding the general case.)

\begin{defn} \label{StandardDefn}
 Suppose 
 \begin{itemize}
 \item $G$ is a subgroup of a Lie group~$H$,
 \item $\Lambda$ is a lattice in~$H$, and
 \item $C$ is a compact subgroup of~$H$ that centralizes~$G$.
 \end{itemize}
 Then $G$ acts (on the right) on the double-coset space $M = \Lambda
\backslash H / C$, and we call~$M$ a \emph{standard $G$-space}.
 \end{defn}

\begin{eg} \label{StandardEg}
 In the setting of Defn.~\ref{StandardDefn}, if $C$~acts freely on
$\Lambda \backslash H$, then $P = \Lambda \backslash H$ is a principal
$C$-bundle over~$M$. In many cases (for example, if $H$ is a
connected, noncompact, simple group), the Mautner phenomenon
\cite{Mautner} implies that $G$~acts irreducibly on~$P$.

There is a $G$-invariant connection on~$P$. To see this, note that,
because $GK \approx G \times K$ is reductive in~$H$, there is an
$\Ad_H(GK)$-invariant complement~$\Lie M$ to~$\Lie C$ in~$\Lie H$. Then
$\Lie M$ defines a $GK$-invariant complement to the vertical tangent
space $\TPvert$.
 \end{eg}
 
We show that if $M$~is a standard
$G$-space, then $P$ must also be a standard $G$-space \see{geomgeneral},
so the possible choices of~$K$ are severely restricted; $K$~must arise
from purely algebraic considerations. 

\begin{thmref}{GeomSS}
 \begin{thm} \label{MainIntro}
 Suppose
 \begin{itemize}
 \item $G$ is a connected, semisimple, linear Lie group, with no compact
factors,
 \item $M = \Lambda \backslash H / C$ is an ergodic, standard $G$-space,
such that $H$~is connected and semisimple, with no compact factors,
 \item $K$ is a connected, compact Lie group,
 and
 \item $P$ is a principal $K$-bundle over~$M$, 
 such that
 \begin{itemize}
 \item the action of~$G$ on~$M$ lifts to an ergodic action of~$G$ on~$P$
by bundle automorphisms,
 and
 \item there is a $G$-invariant connection on~$P$.
 \end{itemize}
 \end{itemize}
 Then there exist
 \begin{itemize}
 \item a finite-index subgroup $\Lambda_0$ of~$\Lambda$,
 \item a connected, compact Lie group~$N$,
 \item a homomorphism $\sigma \colon \Lambda_0 \to N$, with dense image,
 \item a quotient $\overline{C}$ of~$C^\circ$, and
 \item a finite subgroup~$F$ of the center of~$K$,
 \end{itemize}
 such that $K/F \iso \overline{C} \times N$.
 \end{thm}
 \end{thmref}

Combining this with the Margulis Superrigidity Theorem
\ref{MargulisBundle} yields the following result.

\begin{thmref}{SimpleFiniteChoice}
 \begin{cor}
 In the setting of Corollary~\ref{MainIntro}, suppose that $\Rrank H \ge
2$, and that $\Lambda$~is irreducible. Then there is a finite set
$\mathcal{K} = \{K_1,\ldots,K_n\}$ of compact, connected Lie groups,
depending only on $H$, such that $K$~is isomorphic to one of the groups
in~$\mathcal{K}$.
 \end{cor}
 \end{thmref}

\begin{acknowledgements}
 This research was partially supported by grants from the National
Science Foundation. 
 D.W.\ would like to thank the \'Ecole Normale Sup\'erieure de Lyon, the
Newton Institute (Cambridge, U.K.), the University of Chicago, and the
Tata Institute of Fundamental Research (Mumbai, India), for their
hospitality while this research was underway, and while the manuscript
was being edited into a publishable form.
 \end{acknowledgements}

\section{The measurable category}

As is well known, Thms.~\ref{MargulisBundle} and~\ref{ZimmerBundle} can be
restated in the language of Borel cocycles (cf.\ \cite[Prop.~4.2.13,
p.~70]{ZimmerBook} and \cite[Thm.~5.2.5, p.~98]{ZimmerBook}).

\begin{defn}
 Suppose a Lie group~$G$ acts measurably on a Borel space~$M$, with
a quasi-invariant measure~$\mu$, and $H$~is a (second countable) locally
compact group. 
 \begin{itemize}
 \item A Borel measurable function
 $\alpha \colon M \times G \to H$ 
 is a \emph{Borel cocycle} if, for each $g,h \in G$, we have
$\alpha(x,gh) = \alpha(x,g) \, \alpha(xg,h)$ for a.e.\ $x \in
M$.
 \item A Borel cocycle~$\alpha$ is \emph{strict} if the equality holds
for every~$x$, instead of only for almost every~$x$.
 \item If $\alpha \colon M \times G \to H$ is a strict Borel cocycle,
then the \emph{skew-product action $M \atimes H$} of~$G$ is the action
of~$G$ on $M \times H$, given by
 $ (m,h) g = \bigl( mg, h \, \alpha(m,g) \bigr)$.
 \end{itemize}
 Any Borel cocycle is equal, almost everywhere, to a strict
Borel cocycle \cite[Thm.~B.9, p.~200]{ZimmerBook}, so, abusing
terminology, we will speak of the skew-product action $M \atimes H$ even
if $\alpha$~is not strict.
 \end{defn}

\begin{thmref}{MargulisBundle}
 \begin{thm} \label{CocycleG/Gamma}
 For any connected, semisimple, linear Lie group~$G$, with $\Rrank G \ge
2$, there is a corresponding finite set $\mathcal{H} =
\{H_1,\ldots,H_n\}$ of connected, compact Lie groups, such that if
 \begin{itemize}
 \item $H$~is a connected, linear, Lie group, 
 \item $\Gamma$ is an irreducible lattice in~$G$,
 \item $\alpha \colon (\Gamma \backslash G) \times G \to H$ is a Borel
cocycle, 
 and
 \item the skew-product action $(\Gamma \backslash G) \atimes H$ is
ergodic,
 \end{itemize}
 then $H$ is isomorphic to one of the groups in~$\mathcal{H}$.
 \end{thm}
 \end{thmref}

\begin{thmref}{ZimmerBundle}
 \begin{thm}[{Cocycle Superrigidity Theorem}] \label{ZimmerCocycle}
 Let 
 \begin{itemize}
 \item $G$ be a connected, semisimple, linear Lie group, with no compact
factors, such that $\Rrank G \ge 2$,
 \item $M$~be an irreducible, ergodic $G$-space with finite invariant
measure,
 \item $H$ be a connected, semisimple, linear Lie group, with no compact
factors,
 and
 \item $\alpha \colon M \times G \to H$ be a Borel cocycle.
 \end{itemize}
 If the skew-product action $M \atimes H$ is ergodic, then $H$ is locally
isomorphic to a normal subgroup of~$G$. 
 \end{thm}
 \end{thmref}

Recall that if $P$ is a principal $K$-bundle over a manifold~$M$, then
any action of~$G$ on~$P$ by bundle automorphisms yields a Borel cocycle
$\alpha \colon M \times G \to K$, such that $P$ is $G$-equivariantly
measurably isomorphic to the skew product $M \atimes K$ (cf.\
\cite[pp.~66--67]{ZimmerBook}). This suggests the following measurable
version of Prob.~\ref{ProbTop}.

\begin{problem}
 Let 
 \begin{itemize}
 \item $G$ be a connected, semisimple, linear Lie group, with no compact
factors, such that $\Rrank G \ge 2$, 
 \item $M$ be a smooth manifold of finite volume, on which $G$~acts
irreducibly, by volume-preserving diffeomorphisms, and
 \item $K$ be a connected, compact Lie group.
 \end{itemize}
 Describe the Borel cocycles $\alpha \colon M \times G \to K$, such that
the skew-product action $M \atimes K$ is ergodic. In particular, are there
choices of~$G$, $M$, and~$K$ for which no such cocycle exists? (We assume
that every $G$-orbit on~$M$ is a null set, for otherwise
Thm.~\ref{CocycleG/Gamma} applies.)
 \end{problem}

\begin{prop} \label{anyK}
 For any noncompact, linear Lie group~$G$, and any compact
{\upshape(}second countable{\upshape)} group~$K$, there is an
irreducible, ergodic $G$-space~$M$ with finite invariant measure, and a
Borel cocycle $\alpha \colon M \times G \to K$, such that the
skew-product $M \atimes K$ is ergodic {\upshape(}and
irreducible{\upshape)}. 
 \end{prop}

\begin{proof}
 Embed $G$ in some $\SL(n,\real)$.  
For each natural number~$k$, 
	let $H_k = \SL(n+k, \real)$
	and let $\Lambda_k$ be a (cocompact) torsion-free lattice in~$H_k$.  
Let  
 $$M_0   = \frac{H_2}{\Lambda_2} \times \frac{H_3}{\Lambda_3}
\times \frac{H_4}{\Lambda_4} \times \cdots $$
and
 $$ K_0 = \SO(2) \times \SO(3) \times \SO(4) \times \cdots
.$$
 The group~$G$ acts diagonally on~$M_0$.  From vanishing of
matrix coefficients \cite[Thm.~2.2.20, p.~23]{ZimmerBook}, we know that
the diagonal embedding of~$G$ is mixing on each finite subproduct of
$M_0$, so (by approximating arbitrary sets by finite unions of cylinders)
we conclude that $G$ is mixing on the entire product $M_0$. 
In particular, every closed, noncompact subgroup of~$G$ is ergodic
on~$M_0$, so $M_0$ is irreducible.

By construction, $H_k$ contains $G \times \SO(k)$, so there
is a free action of~$K_0$ on~$M_0$ that centralizes the
diagonal action of~$G$.  Then, because $K$ is isomorphic to a
closed subgroup of~$K_0$, we know that there is a free action
of~$K$ on~$M_0$ that centralizes the diagonal action of~$G$. 
Let $M = M_0/K$, so $M_0$ is a principal $K$-bundle over~$M$.
 \end{proof}

\begin{cor} \label{unifK}
 \begin{enumerate}
 \item \label{unifK-mfld}
 If $K$ is a Lie group, then the $G$-space~$M$ can be taken to be a
smooth, compact manifold, and the skew product $M \atimes K$ can be
realized as a principal $K$-bundle over~$M$.
 \item \label{unifK-Lie}
 If $K$ is restricted to the class of Lie groups, then there is a
$G$-space~$M_1$ that works for all $K$ simultaneously.  {\upshape(}That
is, there is no need to vary~$M$ as $K$ ranges over the set of compact
Lie groups.{\upshape)}
 \item \label{unifK-abel}
 Similarly, if $K$ is restricted to the class of compact, abelian groups
that are not necessarily Lie, then there is a $G$-space~$M_2$ that works
for all $K$ simultaneously.
 \end{enumerate}
 \end{cor}

\begin{proof}
 \pref{unifK-mfld} Because $K \subset \SO(k)$, for some~$k$, we may use
$\Lambda_k \backslash H_k$, which is a smooth manifold, instead of~$M_0$,
in the construction.

 \pref{unifK-Lie} Let $K_1$ be the direct
product of one representative from each isomorphism class of
compact Lie groups.  There are only countably many compact
Lie groups, up to isomorphism \cite[Cor.~10.13]{AdamsReduction},
 so $K_1$ is compact and second countable.  From the proposition, there
is an ergodic $G$-space $M_1$ with finite invariant measure and a
$K_1$-valued cocycle~$\alpha_1$ with ergodic skew product.  Now, for
any~$K$, we simply let $\alpha$ be the composition of~$\alpha_1$ with the
projection from $K_1$ to~$K$.

	\pref{unifK-abel} By the argument of~\pref{unifK-Lie}, it suffices to
show that there is a compact abelian group~$A_0$, such that every compact
abelian group is a quotient of~$A_0$.  To see this, let $A$ be the direct
sum of countably many copies of~$\rational$ and countably many copies of
$\rational/\integer$, then let $A_0$ be the Pontryagin dual of~$A$. 
Every countable abelian group is isomorphic to a subgroup of~$A$ (because
every countable abelian group is isomorphic to a subgroup of a countable,
divisible, abelian group \cite[Thm.~24.1, p.~106]{Fuchs} and every
countable, divisible, abelian group is isomorphic to a subgroup of~$A$
\cite[Thm.~23.1, p.~104]{Fuchs}), so, by duality, every compact abelian
group is isomorphic to a quotient of~$A_0$.
 \end{proof}

The following observation is well known, but does not seem to have
previously appeared in print.

\begin{rem}
 Assume $G$~is noncompact and has Kazhdan's property~$T$ (for example,
let $G = \SL(3,\real)$ \cite[Thm.~7.4.2, p.~146]{ZimmerBook}), and let
$\torus = \real/\integer$ be the circle group. From Prop.~\ref{anyK}, we
know there is an ergodic $G$-space~$M$, with finite invariant measure,
and a Borel cocycle $\alpha \colon M \times G \to \torus$, such that $M
\atimes \torus$ is ergodic. Then $\alpha$~is \textbf{not} cohomologous to
any cocycle~$\beta$, such that $\beta(M \times G)$ is countable. To see
this, we argue by contradiction:  let $A$ be the subgroup of~$\torus$
generated by $\alpha(M \times G)$, and suppose, after replacing $\alpha$
with a cohomologous cocycle, that $A$~is countable. Then we may think
of~$\alpha$ as a cocycle into~$A$ (with the discrete topology on~$A$).
Because $G$~has Kazhdan's property~$T$, and $A$~is abelian, we conclude
that $\alpha$ is cohomologous to a cocycle whose values lie in a compact
(thus, finite) subgroup of~$A$ \cite[Thm.~9.1.1, p.~162]{ZimmerBook}.
Thus, we may assume $\alpha(M \times G)$ is finite. But then $M \atimes
\torus$ is clearly not ergodic. This contradicts the choice of~$\alpha$.
 \end{rem}

\section{The geometric category}

 \begin{thm} \label{geometric}
 Let 
 \begin{itemize}
 \item $H$ be a connected Lie group,
 \item $M = \Lambda \backslash H$, for some lattice~$\Lambda$
in~$H$, such that $H$ acts faithfully on~$M$, 
 \item $G$ be a connected, semisimple Lie subgroup of~$H$, with
no compact factors,
 \item $K$ be a compact Lie group,
 \item $E \to M$ be a smooth principal
$K$-bundle, such that the action of~$G$ on~$M$ lifts to a
well-defined {\upshape(}faithful{\upshape)} action of a
cover~$G'$ of~$G$ by bundle automorphisms of~$E$.
 \end{itemize}
 Assume that
 \begin{itemize}
 \item $G$ is ergodic on~$M$, and
 \item $G'$ preserves a connection on~$E$.
 \end{itemize}
 Then there exist
 \begin{itemize}
 \item a closed subgroup~$N$ of~$K$, 
 \item a $G'$-invariant principal $N$-subbundle~$E'$ of~$E$, and
 \item a Lie group~$H'$, with only finitely many connected
components,
 \end{itemize}
 such that 
 \begin{enumerate} 
 \item $H'$ is a transitive group of diffeomorphisms of~$E'$,
 \item $H'$ contains~$N$ as a normal subgroup,
 \item $H'/N$ is isomorphic to~$H$,
 \item \label{geometric-actonM}
 the action induced by $H'$ on $E'/N = M$ is the action of~$H$
on~$M$,
 \item $H'$ contains~$G'$, and
 \item $G'$ is ergodic on~$E'$.
 \end{enumerate}
 Therefore, there is a lattice $\Lambda'$ in~$H'$, such that
 \begin{enumerate} \renewcommand{\theenumi}{\alph{enumi}}
 \item $E'$ is $G'$-equivariantly diffeomorphic to
$\Lambda' \backslash H'$, and
 \item $M$ is $G$-equivariantly diffeomorphic to $\Lambda'
\backslash H' / N$.
 \end{enumerate}
 \end{thm}

\begin{rem}
 Because $K$ centralizes~$G'$, we see that $G'$ is ergodic on $E' k$,
for every $k \in K$. Therefore,
 $\{\, E' k \mid k \in K \,\}$
 is an ergodic decomposition for the $G'$-action on~$E$. In particular,
if $G'$ is ergodic on~$E$, then $E' = E$, so $N = K$.
 \end{rem}

\begin{rem}
 The proof shows that if $K$~is connected and $G'$~is ergodic
on~$E$, then $H'$ may be taken to be connected. In general,
$H'$~is constructed so that $H' = (H')^\circ N$.
 \end{rem}

\begin{proof}
 We begin by establishing some notation. 
 \begin{itemize}
 \item Let $\omega$ be the connection form of some
$G'$-invariant connection \cite[pp.~63--64]{KobayashiNomizu1}.
 \item Let $\Omega = D\omega$ be the curvature form of the
connection \cite[pp.~77]{KobayashiNomizu1}. 
 \item We view $\Lie H$ as the Lie algebra of left-invariant
vector fields on~$H$, so each element of~$\Lie H$ is
well-defined as a vector field on $\Lambda \backslash H = M$. 
 \item For each $X \in \Lie H$, we use $\overline X$ to denote
the lift of~$X$ to a horizontal vector field on~$E$.
 \item For $Z \in \Lie K$, let $\check Z$ be the corresponding
vertical vector field on~$E$ induced by the action of~$K$ (so
$\omega( \check Z ) = Z$).
 \item For any $X \in \Lie G$, let $X'$ be the corresponding
vector field on~$E$ induced by the action of~$G'$.
 \end{itemize}
 By definition, $\Omega$ is a horizontal 2-form on~$E$ taking
values in~$\Lie K$. For each $e \in E$, the connection provides
an identification of the horizontal part of the tangent space
$T_e E$ with~$\Lie H$; therefore, we may think of~$\Omega$ as a
map $\hat\Omega \colon E \to \Hom_{\real}(\Lie H \wedge \Lie H,
\Lie K)$. 
 From \cite[Prop.~II.5.1(c), p.~76 and
comments on p.~77]{KobayashiNomizu1}, we have
 \begin{equation} \label{Kequivariant}
 \mbox{$\hat\Omega_{ek} = (\Ad_K k^{-1}) \circ \hat\Omega_e$
 for every $e \in E$ and $k \in K$.}
 \end{equation}

\setcounter{step}{0}

\begin{step} \label{geometric-orbit}
 We have
 \begin{enumerate}
 \item \label{geometric-orbit-inv}
 $\hat\Omega_{eg} =
\hat\Omega_e$ for every $e \in E$ and $g \in G'$; and
 \item \label{geometric-orbit-inorbit}
 there is some $\phi \in \Hom_{\real}(\Lie H
\wedge \Lie H, \Lie K)$, such that
 $\hat\Omega_e \in \{\, (\Ad_K k) \circ \phi \mid k \in K \,\}$,
 for every $e \in E$.
 \end{enumerate}
 \end{step}
 For any $X \in \Lie H$ and $g,a \in H$, and defining $R_g
\colon H \to H$ by $R_g(h) = hg$, the left-invariance of~$X$
implies
 $d(R_g)_a(X_a) = \bigl( (\Ad_H g^{-1}) X \bigr)_{ga}$
\cite[Prop.~I.5.1, p.~51]{KobayashiNomizu1}; therefore, for $g
\in G'$ and $e \in E$, 
 \begin{equation} \label{horizontalpart}
 \mbox{the horizontal part of $dg_e(\overline X_e)$ is
$\bigl(\overline{(\Ad_H g^{-1}) X}\bigr)_{eg}$.}
 \end{equation}
 Then, because the connection is $G'$-invariant (and $G'$
commutes with the $K$-action on~$E$), we have
 \begin{equation} \label{G'inv}
 \matrix{
 \hat\Omega_e(X,Y)
 &=& \Omega_{e} ( \overline X_e, \overline Y_e) \hfill
 & \mbox{(definition of $\hat\Omega$)} \hfill \cr
 &=& \Omega_{eg} \bigl( dg_e(\overline X_e), dg_e(\overline
Y_e) \bigr) \hfill
 & \mbox{\cite[Prop.~II.6.1(b),
pp.~79--80]{KobayashiNomizu1}} \hfill \cr
 &=& \Omega_{eg} \bigl[ \bigl(\overline{(\Ad_H g^{-1})
X}\bigr)_{eg},
 \bigl(\overline{(\Ad_H g^{-1}) Y}\bigr)_{eg} \bigr] \hfill
 & \mbox{(\eqref{horizontalpart} and $\Omega$ is horizontal)}
\hfill \cr
 &=& \hat\Omega_{eg} \bigl( (\Ad_H g^{-1}) X,
 (\Ad_H g^{-1}) Y \bigr) \hfill 
 & \mbox{(definition of $\hat\Omega$)}  , \hfill 
 }
 \end{equation}
 so $\hat \Omega$ is $G'$-equivariant.
 Then, since the action of~$G'$ on~$E$ preserves a probability measure,
but the action on $\Hom_{\real}(\Lie H \wedge \Lie H, \Lie K)$ is
algebraic, a standard argument based on the Borel Density Theorem
(combining \cite[Thm.~3.1.3, p.~34, Props.~2.1.11 and 2.1.12, p.~11, and
Rem.\ after Thm.~3.2.5, p.~42]{ZimmerBook}) implies that $\hat\Omega$ is
constant on almost every ergodic component of the $G'$-action on~$E$.
 Because almost
every $G$-orbit on~$M$ is dense, this implies that we may choose some
$e_0 \in E$, such that $\hat\Omega$ is constant on $e_0 G'$, and the
projection of $e_0 G'$ to~$M$ is dense.

Let $E'$ be the closure of $e_0 G'$, and let $\phi = \hat\Omega_{e_0}$.
Because $\hat\Omega$ is continuous, we know
that $\hat\Omega$ is constant on~$E'$. Because $K$ is compact,
the projection $\phi \colon E \to M$ is a proper map, so
$\pi(E')$ is closed in~$M$. Since it is also dense, we conclude
that $\pi(E') = M$; thus, $E = E' K$. Now, for any $e \in E'$
and $k \in K$, we have
 $$ \matrix{
 \hat \Omega_{ek} &=& (\Ad_K k^{-1}) \circ \hat\Omega_e \hfill
 & \mbox{(see \eqref{Kequivariant})} \hfill \cr
 &=& (\Ad_K k^{-1}) \circ \hat\Omega_{e_0} \hfill
 & \mbox{($\hat\Omega$ is constant on~$E'$)} \hfill \cr
 &=& (\Ad_K k^{-1}) \circ \phi \hfill
 & \mbox{(definition of~$\phi$).} \hfill
 }
 $$
 This implies that 
 \begin{enumerate}
 \item[\pref{geometric-orbit-inv}] $\hat\Omega$ is constant
on~$E'k$, for each $k \in K$, so $\hat\Omega$ is $G'$-invariant,
and
 \item[\pref{geometric-orbit-inorbit}] $\hat\Omega \in \{\,
(\Ad_K k) \circ \phi \mid k \in K \,\}$,
 for every $e \in E'K = E$.
 \end{enumerate}

\begin{step} \label{omegaconstant}
 We may assume that $\hat\Omega$ is constant.
 \end{step}
 Let 
 \begin{itemize}
 \item $\phi$ be as in Step~\ref{geometric-orbit}, 
 \item $E' = \{\, e \in E \mid \hat\Omega_e = \phi \,\}$,
and
 \item $N = \{\, k \in K \mid (\Ad_K k) \circ \phi = \phi
\,\}$.
 \end{itemize}
 The Implicit Function Theorem (together with
\eqref{Kequivariant} and
Conclusion~\pref{geometric-orbit-inorbit} of
Step~\ref{geometric-orbit}) implies that $E'$ is a smooth
submanifold of~$E$. Then, from \eqref{Kequivariant}, we see
that $E'$ is a principal $N$-subbundle of~$E$. Because
$\hat\Omega$ is $G'$-invariant, this subbundle is
$G'$-invariant.

Because $N$ is compact, we know that $\Ad_K N$ is completely
reducible, so there is an $(\Ad_K N)$-equivariant projection
$p \colon \Lie K \to \Lie N$; let $\omega' = p \circ \omega$.
Then $\omega'$ is a $\Lie N$-valued $1$-form. Because $p$ is
$(\Ad_K N)$-equivariant, we have 
 $R_k^* \omega' = (\Ad_{N} k^{-1}) \circ \omega'$ for all $k
\in N$, so \cite[Prop.~II.1.1, p.~64]{KobayashiNomizu1}
implies that $\omega'$ is the connection form of a connection
on~$E'$. Since the connection corresponding to~$\omega$ is
$G'$-invariant, and $G'$ centralizes~$K$, it is easy to see
that the connection corresponding to~$\omega'$ is also
$G'$-invariant.

If $N \neq K$, then, by replacing $E$ with the subbundle~$E'$, we may
replace $K$ with a smaller subgroup. By the descending chain condition
on closed subgroups of~$K$, this cannot continue indefinitely.

Thus, we may assume $N = K$. Then $E' = E$, so $\hat\Omega$ is
constant on~$E$, as desired.

\begin{step}
 The vector space $\overline{\Lie H} + \check{\Lie K}$ is a Lie
algebra.
 \end{step}
  Note that $[\check {\Lie K}, \overline{\Lie H}] = 0$ (cf.\
\cite[Prop.~II.1.2, p.~65]{KobayashiNomizu1}), so we need
only show $[\overline{\Lie H}, \overline{\Lie H}] \subset
\overline{\Lie H} + \check{\Lie K}$. Fix $X,Y \in \Lie H$,
and let $Z = -2\hat\Omega(X,Y)$. Then, from
\cite[Cor.~II.5.3, p.~78]{KobayashiNomizu1}, we have
 $$ \omega \bigl( [\overline X, \overline Y] \bigr)
 = -2 \Omega \bigl( \overline X, \overline Y \bigr)
 = -2 \hat\Omega(X,Y) 
 = Z
 = \omega( \check Z), $$
 so we see that 
 $$ \bigl[ \overline X, \overline Y \bigr] = \overline{[X,Y]} +
\check Z
 \in \overline{\Lie H} + \check{\Lie K} ,$$
 as desired.

\begin{step}
 We have $\Lie G' \subset \overline{\Lie H}$.
 \end{step}
 For $X,Y,Z \in \Lie G$, and letting $\sum$ denote summation
over the set of cyclic permutations of $(X,Y,Z)$, we have the
following well-known calculation:
 $$ 
 \matrix{
 -\sum \hat\Omega \bigl( [X,Y], Z \bigr) \hfill 
 &=& -\sum \Omega \bigl( \overline{[X, Y]} , \overline Z \bigr)
\hfill
 & \mbox{(definition of~$\hat\Omega$)} \hfill
\cr
 &=& -\sum \Omega \bigl( [\overline X, \overline Y] , \overline
Z \bigr) \hfill
 & \mbox{($\Omega$ is horizontal)} \hfill \cr
 &=& -\sum \Omega \bigl( [\overline X, \overline Y] , \overline
Z \big) \hfill \cr
&& {} + \sum \overline X \bigl( \Omega(\overline Y,\overline
Z) \bigr)\hfill
 & \mbox{($\Omega(\overline Y,\overline
Z) =
\hat\Omega(Y,Z)$ is constant)} \hfill \cr
 &=& 3 \, d\Omega \bigl( \overline X, \overline Y, \overline Z
\bigr) \hfill
 & \mbox{\cite[Prop.~II.3.11, p.~36]{KobayashiNomizu1}}
\hfill \cr
 &=& 3 \, D\Omega \bigl( \overline X, \overline Y, \overline Z
\bigr) \hfill
 & \mbox{($\overline X$, $\overline Y$, and $\overline Z$ are
horizontal)} \hfill \cr
 &=& 0 \hfill & \mbox{(Bianchi
Identity \cite[Thm.~II.5.4, p.~78]{KobayashiNomizu1}).} \hfill
 } $$
 Thus, as is well known, $\hat\Omega|_{\Lie G \wedge \Lie G}$ is
a 2-cocycle for the Lie algebra cohomology of~$\Lie G$ (with
coefficients in the module~$\Lie K$ with trivial $\Lie
G$-action) \cite[(3.12.5), p.~220]{Varadarajan}. From
Whitehead's Lemma \cite[Thm.~3.12.1, p.~220]{Varadarajan}, we
know that $H^2(\Lie G; \Lie K) = 0$, so there is some $\sigma
\in \Lie \Hom_{\real}(\Lie G, \Lie K)$, such that
 $$ \mbox{$\hat\Omega(X,Y) = \sigma \bigl( [X,Y] \bigr)$ for all
$X,Y \in \Lie G$.} $$

For $g \in G$ and $X,Y \in \Lie G$, we have
 $$
 \matrix{
 \sigma \bigl( [X,Y] \bigr)
 &=& \hat\Omega(X,Y) \hfill
 & \mbox{(definition of~$\sigma$)} \hfill  \cr
 &=& \hat\Omega\bigl( (\Ad_H g^{-1})X, (\Ad_H g^{-1})Y \bigr) \hfill
 & \mbox{(\eqref{G'inv} and $\hat\Omega$ is constant)} \hfill \cr
 &=& \sigma \bigl( [(\Ad_H g^{-1})X, (\Ad_H g^{-1})Y] \bigr) \hfill
 & \mbox{(definition of~$\sigma$)} \hfill \cr
 &=& \sigma \bigl( (\Ad_H g^{-1})[X, Y] \bigr) 
 & \mbox{($\Ad_H g^{-1}$ is an automorphism of $\Lie H$),}
\hfill
 }
 $$
  Then, because
$[\Lie G, \Lie G] = \Lie G$, we conclude that
 $ \sigma \bigl( (\Ad g)X - X \bigr) = 0$, 
 for every $g \in G$ and $X \in \Lie G$. Because $G$ is
semisimple, we have $[G,\Lie G] = \Lie G$, so we conclude
that $\sigma(\Lie G) = 0$. 
 Combining this with the definitions of~$\sigma$
and~$\hat\Omega$, and the fact that the form~$\Omega$ is
horizontal, we have
 \begin{equation} \label{Omega(G,G)=0}
  0 = \sigma \bigl( [X,Y] \bigr)
 = \hat\Omega(X,Y)
 = \Omega \bigl( \overline X, \overline Y \bigr)
 = \Omega(X', Y') ,
 \end{equation}
 for every $X,Y \in \Lie G$. 

For any fixed $e_0 \in E$, define $\tau \colon \Lie G \to
\Lie K$ by $\tau(X) = \omega \bigl( X' (e_0) \bigr)$. 
 Then, because 
 $$2\Omega(X', Y') = [\tau(X), \tau(Y)] - \tau \bigl( [X,Y]
\bigr)$$
 \cite[Prop.~II.11.4, p.~106]{KobayashiNomizu1}, we conclude,
from \pref{Omega(G,G)=0}, that $\tau \colon \Lie G \to \Lie K$
is a homomorphism. Since $G$ has no compact factors, this
implies $\tau = 0$. Because $e_0$ is an arbitrary element
of~$E$, we conclude that $\omega(\Lie G') = 0$,  so every
vector field in $\Lie G'$ is horizontal. Therefore $\Lie G' =
\overline{\Lie G} \subset \overline{\Lie H}$.

\begin{step} \label{proofcompletion}
 Completion of the proof.
 \end{step}
 Let $H'_0$ be the connected Lie group of diffeomorphisms of~$E$
corresponding to the Lie algebra $\overline{\Lie H} +
\check{\Lie K}$. Because $K$ centralizes~$\overline{\Lie H}$
\cite[Prop.~II.1.2, p.~65]{KobayashiNomizu1} and (obviously)
$K$ is normalized by~$K$, we see that $K$ is normalized by
$H'_0$. Let $H' = H'_0 K$ (so $H'_0$ is the identity component
of $H'$, and $H'_0$ is a finite-index subgroup of~$H'$); then
$K$ is a normal subgroup of~$H'$. Furthermore, because $\Lie G'
\subset \overline{\Lie H} \subset \Lie H'$, we have $G' \subset
H'$.

 Because $\Lie H' = \Lie H'_0 = \overline{\Lie H} + \check{\Lie
K}$, we see that $H'$ is transitive on~$E$, with discrete
stabilizer, so there is a discrete subgroup~$\Lambda'$ of~$H'$,
such that $E$ is $H'$-equivariantly diffeomorphic to $\Lambda'
\backslash H'$. Then, because $G' \subset H'$, we know that $E$
is $G'$-equivariantly diffeomorphic to $\Lambda' \backslash H'$.
Then, because 
 $$ \Lie H'/ \check{\Lie K}
 = (\overline{\Lie H} + \check{\Lie K})/ \check{\Lie K}
 \iso \Lie H $$
 and $H' = H'_0 K$, we see that the action induced by~$H'$ on
$\Lambda' \backslash H' / K \iso E/K = M$
 is the action of~$H$ on~$M$. This implies that $H'/K$ is a
cover of~$H$. Also, because there is an $H$-invariant
probability measure on $M = \Lambda \backslash H$, and $K$ is
compact, this implies that there is an $H'$-invariant
probability measure on $\Lambda' \backslash H'$; thus,
$\Lambda'$ is a lattice in~$H'$.

Let $L$ be the smallest normal subgroup of $H'_0$ that
contains~$G'$. (Note that $L$ is also normal in~$H'$ (because
$H' = H'_0 K$ and $K$ centralizes~$G'$), so $\Lambda' L$
is a subgroup of~$H'$. Let
 \begin{itemize}
 \item $H''$ be the identity component of the closure of
$\Lambda' L$,
 \item $N' = K \cap \Lambda' H''$, and
 \item $E'' = \Lambda' \backslash \Lambda' H'' \subset E$.
 \end{itemize}
 Note that $\Lambda'$ normalizes $H''$, so $\Lambda' H''$ is a
(closed) subgroup of~$H'$.

Because $G$ is ergodic on~$M$, we may assume that $\Lambda' G'
K$ is dense in~$H'$. Then, because $\Lambda' H''$ is closed and
contains $\Lambda' G'$, we conclude that $\Lambda' H'' K = H'$.
This means $E'' K = E$, so $E''$ projects onto all of~$M$.

We claim that $E''$ is a principal $N'$-bundle over~$M$. It is
clear that $E''$ is $N'$-invariant, so $E''$ is a principal
$N'$-bundle over $E''/N'$. Thus, we need only show that the
natural map $E''/N' \to E/K = M$ is injective. For $\lambda_1,
\lambda_2 \in \Lambda'$, $h_1,h_2 \in H''$, and $k \in K$, such
that $\lambda_1 h_1 k = \lambda_2 h_2$, we have
 $$ k = (\lambda_1 h_1)^{-1} (\lambda_2 h_2) \in \Lambda' H''
,$$
 so $k \in N'$, as desired.

Arguing as in the last three paragraphs of
Step~\ref{omegaconstant}, we see that we may assume $N' = K$.
Then $E'' = E$, so $\Lambda' L$ is dense in $\Lambda'
\backslash H'$. From the Mautner Phenomenon \cite{Mautner}, we
conclude that $G'$ is ergodic on $\Lambda' \backslash H' = E$.

All that remains is to show that $H'/K$ acts faithfully on~$M$.
(This implies that $H'/K$ is isomorphic to~$H$.) Thus, we wish
to show that $K$ is the kernel of the action of~$H'$ on
$\Lambda' \backslash H' / K$. Now $H'$ is transitive on
$\Lambda' \backslash H' / K$, with stabilizer $\Lambda' K$.
Thus, the kernel of the action is a (normal) subgroup of the
form $LK$, where $L$ is a subgroup of~$\Lambda'$. (We wish to
show $L$ is trivial.) Let $\lambda \in L$. Then, because $L$ is
discrete and $G'$ is connected, we must have $\lambda^{G'}
\subset \lambda K$. Since $K$ is compact, but $G'$ has no
compact factors, we conclude that $\lambda^{G'} = \lambda$.
Therefore $L^{G'} = L$, so
 $L^{G' \Lambda'} = L^{\Lambda'} \subset \Lambda'$.
 Because $G'$ is ergodic on $E = \Lambda' \backslash H'$, we
may assume that $G' \Lambda'$ is dense in~$H'$. So $L^{H'}
\subset \Lambda'$. Since $H'$ is faithful on $\Lambda'
\backslash H'$, we know $\Lambda'$ does not contain any
nontrivial normal subgroups of~$H'$. Therefore $L^{H'}$ must be
trivial, so $L$ is trivial, as desired.
 \end{proof}

\begin{rem}
 If $G'$ is ergodic on~$E$, and the center of~$K$ is discrete,
then combining Step~\ref{omegaconstant} with
\eqref{Kequivariant} yields the conclusion that $\hat\Omega =
0$; thus, under these assumptions, every $G'$-invariant
connection on~$E$ is flat.
 \end{rem}

The following corollary generalizes Thm.~\ref{geometric}, by allowing $M$
to be a double-coset space $\Lambda \backslash H / C$, instead of
requiring it to be a homogeneous space $\Lambda \backslash H$. However,
we add the additional assumption that $G$~is ergodic on the bundle~$P$.
See Rem.~\ref{techassumps} for a discussion of some other less important
assumptions in the statement of this more general result.

 \begin{cor} \label{geomgeneral}
 Let 
 \begin{itemize}
 \item $H$ be a connected Lie group,
 \item $\Lambda$ be a torsion-free lattice in~$H$, such that $H$
acts faithfully on~$\Lambda \backslash H$, 
 \item $C$ be a compact subgroup of~$H$,
 \item $M = \Lambda \backslash H / C$,
 \item $G$ be a connected, semisimple Lie subgroup of~$H$, with
no compact factors,
 \item $K$ be a connected, compact Lie group,
 \item $P \to M$ be a smooth principal
$K$-bundle, such that the action of~$G$ on~$M$ lifts to a
well-defined {\upshape(}faithful{\upshape)} action of a
cover~$G'$ of~$G$ by bundle automorphisms of~$P$.
 \end{itemize}
 Assume that
 \begin{itemize}
 \item $G$ is ergodic on $\Lambda \backslash H$,
 \item $G$ is ergodic on~$P$, and
 \item $G'$ preserves a connection on~$P$.
 \end{itemize}
 Then there exist
 \begin{itemize}
 \item a Lie group~$H'$, with only finitely many connected
components,
 \item a lattice $\Lambda'$ in~$H'$,
 \item a compact subgroup $C'$ of~$H'$, 
 \item a $G'$-equivariant diffeomorphism $\phi \colon \Lambda'
\backslash H' / C' \to P$,
 and
 \item a continuous, surjective homomorphism $\rho \colon H' \to
H$, with compact kernel,
 \end{itemize}
 such that, letting
 $$ K' = \rho^{-1}(C) ,$$
 we have
 \begin{enumerate} 
 \item $H'$ contains $G'$, and $G'$~acts ergodically on
$\Lambda' \backslash H'$,
 \item $G'$ centralizes~$K'$, and $C'$~is a normal subgroup
of~$K'$,
 \item $\phi$ conjugates the action of $K'/C'$ on $\Lambda'
\backslash H' / C'$ to the action of~$K$ on~$P$,
 and
 \item $\phi$ factors through to a diffeomorphism that
conjugates the action of~$G'$ on $\Lambda' \backslash H' / K'$
to the action of~$G$ on $M = \Lambda \backslash H /K$.
 \end{enumerate}
 \end{cor}

\begin{proof}
 Let
 \begin{itemize}
 \item $\sigma \colon \Lambda \backslash H \to \lambda
\backslash H / C$ be the natural quotient map;
 \item $\zeta \colon P \to \Lambda \backslash H / C$ be the
bundle map; and
 \item $E = \{\, (x,p) \in (\Lambda \backslash H) \times P \mid
\sigma(x) = \zeta(p) \,\}$.
 \end{itemize}
 Then $E$ is the principal $K$-bundle over~$\Lambda \backslash
H$ obtained as the pullback of~$P$. Note that
 \begin{itemize}
 \item $G'$ acts (diagonally) on~$E$, via $(x,p)g = (x g, p g)$,
 and
 \item $K$ acts on~$E$, via $(x,p)k = (x,pk)$.
 \end{itemize}
 \emph{Warning:} $G'$ may not be ergodic on~$E$ \see{G'noterg}.

\setcounter{step}{0}

\begin{step}
 $G'$ preserves a connection on~$E$.
 \end{step}
 Let $\omega$ be the connection form of a $G'$-invariant
connection on~$P$, and define
 $\omega^E_{(x,p)}(v,w) = \omega_p(w) \in \Lie K$ 
 for $e = (x,p) \in E$ and 
 $$(v,w) \in T_{(x,p)}(E) \subset T_x(\Lambda \backslash H)
\oplus T_p(P) .$$
 Because $\omega$ is a connection form, it is easy to see, from
the characterization of connection forms \cite[Prop.~2.1.1,
p.~64]{KobayashiNomizu1}, that $\omega^E$ is the connection
form of a connection on~$E$.
 Because $\omega$ is $G'$-invariant, it is clear that
$\omega^E$ is $G'$-equivariant.

\begin{step}
 Let 
 \begin{itemize}
 \item $N$, $E'$, $H'$, and~$\Lambda'$ be as in
Theorem~\ref{geometric}, 
 \item $\psi \colon \Lambda' \backslash H' \to E'$ be
an $H'$-equivariant diffeomorphism, so $\psi \colon \Lambda'
\backslash H' \hookrightarrow E$,
 \item $\rho \colon H' \to H$ be the surjective homomorphism,
with kernel~$N$, that results from \fullref{geometric}{actonM},
 and
 \item $\zeta' \colon \Lambda' \backslash H' \to \Lambda
\backslash H$ be the affine map induced by~$\rho$.
 \end{itemize}
 \end{step}

\begin{equation} \label{geomgeneral-commdiag}
 \begin{CD} 
 H' @>>> \Lambda' \backslash H' @>\psi>> E @>{\pi_2}>> P \\
 @VV\rho V @VV{\zeta'}V @VV{\pi_1}V @VV\zeta V \\
 H @>>> \Lambda \backslash H @= \Lambda \backslash H @>\sigma>>
\Lambda \backslash H / C @= M \\
 \end{CD} 
 \end{equation}

\begin{step}
 $K' = \rho^{-1}(C)$ is compact, and $G'$ centralizes~$K'$.
 \end{step} 
 Because $G$ normalizes~$C$, we know, by pulling back to~$H'$,
that $G'$ normalizes~$K'$. Furthermore, because both $N$ and
$K'/N \iso C$ are compact, we know that $K'$ is compact. Since
$G'$ has no compact factors, and normalizes~$K'$, we conclude
that $G'$ centralizes~$K'$.

\begin{step}
 $\psi$ factors through to a diffeomorphism that conjugates the
action of~$G'$ on $\Lambda' \backslash H' / K'$ to the action
of~$G$ on $\Lambda \backslash H / C$. 
 \end{step}
 From the definition of~$\rho$, we see that the
$G'$-equivariant map $\psi$ induces a diffeomorphism that
conjugates the action of~$K'$ on $\Lambda' \backslash H' / N$
to the action of $\rho(K') = C$ on $\Lambda \backslash H$. Thus,
it factors through to a $G'$-equivariant diffeomorphism between
the orbit spaces of these actions.

\begin{step} \label{geomgeneralpf-C0}
 There is a closed subgroup~$C'$ of~$K'$, such that for
$e_1,e_2 \in E'$, we have $\pi_2(e_1) = \pi_2(e_2)$ if and only
if $e_1 \in e_2 C'$; hence $\pi_2 \circ \psi$ factors through
to a $G'$-equivariant diffeomorphism $\phi \colon \Lambda'
\backslash H' / C' \to P$.
 \end{step}
 For $e \in E'$, let
 $$ \phi(e) = \{\, c \in K' \mid \pi_2(e c) = \pi_2(e) \,\} .$$
 Because $G'$ centralizes~$K'$, and because, from its
definition, $\pi_2$ is obviously $G'$-equivariant, we have
$\phi(e g) = \phi(e)$ for $e \in E'$ and $g \in G'$. Then,
because $G'$ is ergodic on~$E'$, we conclude that $\phi$ is
essentially constant. By continuity, it must be constant: let
$C' = \phi(e)$ for any $e \in E'$.

From the definition of~$\phi$, we see that if $e_1 \in e_2
C'$, then $\pi_2(e_1) = \pi_2(e_2)$. Conversely, suppose
$\pi_2(e_1) = \pi_2(e_2)$. Write $e_1 = \psi(\Lambda' h_1)$ and
$e_2 = \psi(\Lambda' h_2)$, for some $h_1, h_2 \in H'$. Because
$\zeta \bigl( \pi_2(e_1) \bigr) = \zeta \bigl( \pi_2(e_3)
\bigr)$, we see, from the commutative diagram
\pref{geomgeneral-commdiag}, that
 $\zeta'(\Lambda' h_1) \in \zeta'(\Lambda' h_2) C$.
 Thus, from the definitions of~$\zeta'$ and~$K'$, we conclude
that 
 $$\Lambda' h_1 \in \Lambda' h_2 \cdot \rho^{-1}(C) =  \Lambda'
h_2 K' .$$
 So
 $$ e_1 = \psi(\Lambda' h_1) \in \psi(\Lambda' h_2 K') =
\psi(\Lambda' h_2) K' = e_2 K' .$$
 Then, because $\pi_2(e_1) = \pi_2(e_2)$, we conclude, from the
definition of~$C'$, that $e_1 \in e_2 C'$.

For $c_1, c_2 \in C'$, and any $e \in E'$, we have
 $$ \pi_2(e c_1 c_2) = \pi_2(e c_1) = \pi_2(e) ,$$
 so $c_1 c_2 \in C'$. Therefore, $C'$ is a subgroup of~$K'$.
Because $H'$ acts continuously, and $\pi_2$ is continuous, we
know that $C'$ is closed.

Because $\pi_2 \circ \psi$ is $G'$-equivariant, and $G'$~is
ergodic on~$P$, we see that $\pi_2 \circ \psi$ is surjective.
Hence
 $$ P = \pi_2(E') \approx E'/C' \approx \Lambda' \backslash H'
/ C' .$$

\begin{step} 
 $C'$ is normal in~$K'$, with $K'/C' \iso K$, and
$\phi$~conjugates the action of $K'/C'$ on $\Lambda' \backslash
H' / C'$ to the action of~$K$ on~$P$.
 \end{step}
 Fix some $e \in E'$. From the commutative diagram
\pref{geomgeneral-commdiag}, and the definition of~$K'$, we see
that $\pi_2(e K')$ is a fiber of~$\zeta$; hence,
Step~\ref{geomgeneralpf-C0} implies that $C' \backslash K'$ is
diffeomorphic to a fiber of~$\zeta$.
 Because $K$ is connected, we know that the fiber is connected,
so we conclude that $C' \backslash K'$ is connected. Therefore
$K' = (K')^\circ C'$. 

Let $X \in \Lie K'$, and let $X'$ be the corresponding vector
field on $\Lambda' \backslash H'$ induced by the action of~$K'$
on $\Lambda' \backslash H'$. Any element~$Z$ of~$\Lie K$
induces a vertical vector field~$\check Z$ on~$P$, and
$\check{\Lie K}$ constitutes the entire vertical tangent space
at each point. Thus, for any fixed $x_0 \in \Lambda' \backslash
H'$, there is some $Z \in \Lie K$, such that we have $d (\pi_2
\circ \psi)_{x_0} (X') = \check Z_{\pi_2(\psi({x_0}))}$. Now,
each of the maps $x \mapsto d (\pi_2 \circ \psi)_x (X')$ and $x
\mapsto \check Z_{\pi_2(\psi({x_0}))}$ is $G'$-equivariant
(because $G'$ centralizes both $K'$ and~$K$), so these maps are
equal on the orbit $x_0 G'$. We may choose $x_0$ so that this
orbit is dense, and then we conclude, by continuity, that 
 \begin{equation} \label{geomgeneralpf-dphi}
 \mbox{$d (\pi_2 \circ \psi)_{x} (X') = \check
Z_{\pi_2(\psi({x_0}))}$ for all $x \in \Lambda' \backslash H'$.}
 \end{equation}
 This implies that $X'$ factors through to a well-defined
vector field on $\Lambda' \backslash H' / C'$. Hence
 $(\Ad c)X \in X + \Lie C'$, for all $c \in C'$. Because $X$~is
an arbitrary element of~$\Lie K'$, this implies that
$(K')^\circ$ normalizes~$C'$. Because $K' = (K')^\circ C'$, we
conclude that $K'$ normalizes~$C'$. Furthermore, because
$d(\pi_2 \circ \psi)$ maps $\Lie K'$ to~$\Lie K$
\see{geomgeneralpf-dphi}, we know that $\pi_2 \circ \psi$
conjugates the action of $K'$ on $\Lambda' \backslash H'$ to
the action of~$K$ on~$P$. Hence, $\phi$ conjugates the action of
$K'/C'$ on $\Lambda' \backslash H' / C'$ to the action of~$K$
on~$P$.
 \end{proof}

\begin{rem} \label{techassumps}
 Some of the technical assumptions in
Corollary~\ref{geomgeneral} are not very important; they can be
satisfied by passing to a finite cover, or making other similar
minor adjustments.
 \begin{enumerate}

 \item One might assume only that $M$ is connected, rather than
that $H$ is connected. Then $\Lambda_0 \backslash H^\circ /
C_0$ is a finite cover of~$M$, where $\Lambda_0 = \Lambda \cap
H^\circ$ and $C_0 = C \cap H^\circ$.

 \item The assumption that $H$ acts faithfully on $\Lambda
\backslash H$ is only a convenience; one could always mod out
the kernel of this action.

 \item If one does not assume that $\Lambda$ is torsion free,
then Selberg's Lemma \cite[Cor.~6.13, p.~95]{RaghunathanBook} implies that
$\Lambda$ has a torsion-free subgroup $\Lambda_0$ of finite
index. The space $\Lambda_0 \backslash H / C$ is a finite cover
of~$M$.

 \item If one does not assume that $K$ is connected, then
$P/K^\circ$ is a finite cover of~$M$.

\item
 If one does not assume that $G$~is ergodic on $\Lambda
\backslash H$, then we can construct a closed, connected
subgroup~$H_0$ of~$H$, such that, after replacing $\Lambda$ by
a conjugate subgroup,
 \begin{itemize}
 \item $H_0$ contains~$G$, \item $(\Lambda \cap H_0) \backslash
H_0 / (C \cap H_0)$ is a finite cover of~$M$, and
 \item $G$ is ergodic on $(\Lambda \cap H_0) \backslash H_0$.
 \end{itemize}
 To see this, write $H = LKR$, where $R$~is the
radical of~$H$, $K$~is the maximal compact, semisimple
quotient of~$H$, and $L$~is a connected, semisimple subgroup
of~$H$, with no compact factors. Let 
 \begin{itemize}
 \item $H_0$ be the identity component of the closure of
$\Lambda L [R,L]$,
 \item $\Lambda_0 = \Lambda \cap H_0$, and
 \item $C_0 = C \cap H_0$.
 \end{itemize}
 We know that $G$ is ergodic on~$\Lambda \backslash H / C$
(because $G'$ is ergodic on~$P$), so, by replacing $\Lambda$
with a conjugate subgroup, we may assume that $\Lambda G C$ is
dense in~$H$; then $\Lambda H_0 C = H$. Let
 \begin{itemize}
 \item $\Lambda_1 = \Lambda \cap (H_0 C)$,
 \item $H_1 = \Lambda_1 H_0$, and
 \item $C_1 = C \cap H_1$. 
 \end{itemize}
 Let us show that the natural map $\Lambda_1 \backslash H_1 /
C_1 \to \Lambda \backslash H / C = M$ is a bijection.
 \begin{itemize}
 \item Because $\Lambda H_0 C = H$, and $H_0 \subset H_1$, we
know that the map is surjective.

 \item Suppose $\lambda h_1 c \in H_1$, with $\lambda \in
\Lambda$, $h \in H_1$, and $c\in C$. Then
 $$ \lambda \in H_1 C H_1
 \subset H_0 C ,$$
 so $\lambda \in \Lambda \cap (H_0 C) = \Lambda_1 \subset H_1$.
Therefore $c \in H_1$. So we conclude that the map is injective.
 \end{itemize}
 Hence, $\Lambda_1 \backslash H_1 / C_1$ is diffeomorphic
to~$M$.

Now $H_1/H_0$ is discrete, and contained in the compact group
$C H_0 / H_0$, so it must be finite. That is, $H_0$ has finite
index in~$H_1$. Therefore, $\Lambda_0$ has finite index
in~$\Lambda_1$, and $C_0$ has finite index in~$C_1$. Hence,
$\Lambda_0 \backslash H_0 / C_0$ is a finite cover of $\Lambda_1
\backslash H_1 / C_1 \iso M$. The Mautner phenomenon
\cite{Mautner} implies that $G$ is ergodic on $\Lambda_0
\backslash H_0$.

 \end{enumerate}
 \end{rem}

\begin{cor} \label{GeomSS}
 In the setting of Corollary~\ref{geomgeneral}, suppose $H$~is
semisimple, with no compact factors. Then there exist
 \begin{itemize}
 \item a finite-index subgroup $\Lambda_0$ of~$\Lambda$,
 \item a connected, compact Lie group~$N$,
 \item a homomorphism $\sigma \colon \Lambda_0 \to N$, with
dense image,
 \item a quotient $\overline{C}$ of~$C^\circ$, and
 \item a finite subgroup~$F$ of the center of~$K$,
 \end{itemize}
 such that $K/F \iso \overline{C} \times N$.
 \end{cor}

\begin{proof}
 Let $N = \ker(\rho)^\circ$, so $(H')^\circ = HN$. Because
$N$~is a compact, normal subgroup of~$H'$, and $H$ has no
compact factors, we see that $(H')^\circ$ is isogenous to~$H
\times N$. By modding out a finite group, let us assume
$(H')^\circ = H \times N$.

Let $\Lambda_0 = \Lambda \cap (H')^\circ$, and let $\sigma
\colon \Lambda_0 \to N$ be the projection into the second
factor of $H \times N$. Because $G'$ has no compact factors, we
must have $G' \subset H$. Since $G'$~is ergodic on $\Lambda'
\backslash H'$, and, hence, on $\Lambda_0 \backslash
(H')^\circ$, this implies that $H \Lambda_0$ is dense in
$(H')^\circ$, so $\sigma(\Lambda_0)$ is dense in~$N$.

We have $(K')^\circ = \rho^{-1}(C)^\circ = C^\circ \times N$,
and $C'_0 = C' \cap (K')^\circ$~is a normal subgroup
of~$(K')^\circ$, such that $(K')^\circ/C'_0 \iso K$. Therefore,
 $$ K \iso \frac{(K')^\circ}{C'_0}
 \approx \frac{C^\circ N}{C'_0 N} \times \frac{N}{N \cap C'_0}
,$$
 where ``$\approx$" denotes isogeny, that is, an isomorphism
modulo appropriate finite groups.
 \end{proof}

\begin{cor} \label{SimpleFiniteChoice}
 In the setting of Corollary~\ref{geomgeneral}, suppose that $H$~is
semisimple, with no compact factors, that $\Rrank H \ge 2$, and that
the lattice~$\Lambda$ is irreducible. Then there is a finite list
$K_1,\ldots,K_n$ of compact groups, depending only on $H$, such that
$K$~is isomorphic to~$K_j$, for some~$j$.
 \end{cor}

\begin{proof}
 We apply Cor.~\ref{GeomSS}. There are only finitely many connected,
compact Lie groups of any given dimension, and $\dim K = (\dim
\overline{C}) + (\dim N)$, so it suffices to find bounds on $\dim
\overline{C}$ and $\dim N$ that depend only on~$H$. We have $\dim
\overline{C} \le \dim H$.  Given $H$, the Margulis Superrigidity Theorem
\pref{MargulisBundle} implies that there are only finitely many choices
for~$N$, up to isomorphism, so $\dim N$ is bounded.
 \end{proof}

In Cor.~\ref{SimpleFiniteChoice}, one may replace the
assumption that $H$~is simple with the weaker assumption that 
 \begin{itemize}
 \item $H$~is semisimple,
 and
 \item the lattice $\Lambda$ is irreducible.
 \end{itemize}

\begin{eg} \label{G'noterg}
 If $M = \Lambda \backslash H / K$ and $P = \Lambda \backslash
H$, then the principal bundle~$E$ constructed in the proof of
Cor.~\ref{geomgeneral} is $G'$-equivariantly diffeomorphic to
$K \times P$, with $(k,p)g = (k,pg)$. So $G'$~is not ergodic
on~$E$ (unless $K$~is trivial).
 \end{eg}

\end{article}


\begin{thebibliography}{9}

\bibitem{AdamsReduction}
 S.~Adams:
 Reduction of cocycles with hyperbolic targets,
 \emph{Ergodic Th. Dynam. Sys.} 16 (1996), no. 6, 1111--1145. 

\bibitem{Fuchs}
 L.~Fuchs:
 \emph{Infinite Abelian Groups, vol.~I},
 Academic Press, New York, 1970.

\bibitem{KobayashiNomizu1}
 S.~Kobayashi and K.~Nomizu:
 \emph{Foundations of Differential Geometry, vol.~1},
 Interscience, New York, 1963.

\bibitem{MargulisBook}
 G.A.~Margulis,
 \emph{Discrete Subgroups of Semisimple Lie Groups,}
 Springer, New York, 1991.

\bibitem{Mautner}
 Moore, C.C.: 
 The Mautner phenomenon for general unitary representations, 
 \emph{Pacific J. Math.} 86 (1980), no. 1, 155--169. 

\bibitem{RaghunathanBook}
 M.S.~Raghunathan:
 \emph{Discrete Subgroups of Lie Groups},
 Springer, New York, 1972.

\bibitem{Varadarajan}
 V.~S.~Varadarajan:
 \emph{Lie Groups, Lie Algebras, and Their
Representations}, Springer-Verlag, New York, 1984.

\bibitem{ZimmerBook}
 R.~J.~Zimmer,
 \emph{Ergodic Theory and Semisimple Groups,}
 Birkh\"auser, Boston, 1984.


\end{thebibliography}
\end{document}